\documentclass{article}

\usepackage[english]{babel}

\usepackage{amsmath}
\usepackage{amssymb}
\usepackage{graphicx}
\usepackage[colorlinks=true, allcolors=blue]{hyperref}

\usepackage{pgfplots}
\usepgfplotslibrary{groupplots}
\usetikzlibrary{spy}
\usetikzlibrary{backgrounds}

\usepackage{multirow}

\pgfplotsset{
    compat=newest,
    height=0.4\textwidth, width=\textwidth,
    tick label style={font=\scriptsize},
	label style={font=\small},
    every axis legend/.append style={font=\scriptsize},
            }
            
\usepackage{listings}
\usepackage[ruled, noline, linesnumbered]{algorithm2e}

\newcommand{\U}{\mathbf{U}}
\newcommand{\F}{\mathbf{F}}

\newcommand{\eps}{\varepsilon}

\renewcommand{\P}{\mathbb{P}}
\newcommand{\softplus}{\text{softplus}}	
\newcommand{\llbracket}{[\![}
\newcommand{\rrbracket}{]\!]}
	
\newcommand{\correc}[2]{\textcolor{red}{#1} \textcolor{blue}{#2}}
\usepackage{comment}
\usepackage[letterpaper, top=2cm, bottom=2cm, left=3cm, right=3cm, marginparwidth=1.75cm]{geometry}
\usepackage{layouts}
\usepackage{authblk}
\title{An optimal control deep learning method to design artificial viscosities for Discontinuous Galerkin schemes}
\author[1, 2]{Léo Bois}
\author[1, 2]{Emmanuel Franck}
\author[1, 2]{Laurent Navoret}
\author[1]{Vincent Vigon}
\affil[1]{Institut de Recherche Mathématique Avancée, UMR 7501 Université de Strasbourg et CNRS, 7 rue René Descartes 67000 Strasbourg, France}
\affil[2]{INRIA Nancy-Grand Est, TONUS Project, Strasbourg, France}

\begin{document}
	\maketitle

	\begin{abstract}
		In this paper, we propose a method for constructing a neural network viscosity in order to reduce the non-physical oscillations generated by high-order Discontiuous Galerkin (DG) methods. To this end, the problem is reformulated as an optimal control problem for which the control is the viscosity function and the cost function involves comparison with a reference solution after several compositions of the scheme. The learning process is strongly based on gradient backpropagation tools. Numerical simulations show that the artificial viscosities constructed in this way are just as good or better than those used in the literature.
    \end{abstract}
	
	\section{Introduction}

In computational fluid dynamics, Discontinuous Galerkin (DG) methods can be used as an alternative to finite volumes (FV) methods when a high order of convergence is desired. Indeed, by using polynomials coupled with the weak form of the equation to approximate the solution, DG methods allow to reach arbitrarily high orders of convergence \cite{cockburn2001runge,hesthaven2007nodal}. However, when the solution exhibits shocks or strong gradients, these high order methods introduce non-physical oscillations, which can deteriorate the accuracy of the solution and lead to stability issues. Since shocks and strong gradients easily appear in non-linear conservative systems, even from continuous initial conditions, countermeasures are required to stabilize DG methods. 

Classical approaches to reduce oscillations and stabilize DG methods are based on slope limiters, filtering techniques or artificial viscosity methods. Slope limiter methods were initially developed for FV methods and then adapted to DG schemes: they use a troubled-cell indicator to identify cells with oscillations and then define fluxes at the cells interfaces with second order polynomial approximation and total variations diminishing property \cite{cockburn1989tvb,cockburn1989tvb_onedim}. An alternative is to consider WENO (Weighted Essentially Non Oscillating) type reconstruction to take advantage of the full high-order approximation in the identified troubled cells and their neighbours \cite{zhong2013simple, qiu2004hermite,zhao2020hybrid}. Regarding filtering techniques, they involve applying a linear filter using the modal representation of the solution locally in each cell to smooth out the solution \cite{hesthaven2008filtering,hesthaven2017numerical}. Finally the artificial viscosity method consists in adding a non-linear viscous term to the equation, which makes the solution smoother, thus saving the DG methods from situations that they cannot handle correctly. The viscous term can then be tuned with a local coefficient, which should only be activated in problematic areas and should vanish as the characteristic length of the mesh tends to zero. 

In this paper, we will focus on artificial viscosity methods. A few different approaches have emerged to deal with different problems. For instance, the artificial viscosity can be function of the divergence of the velocity field \cite{mani2009suitability}, function of the modal decay of the solution in each cell \cite{persson2006sub} or function of the entropy production \cite{guermond2011entropy}. A comparative study of some models has been carried out in \cite{Yu:276634}, and shows that these different models behave differently from each other, and that the most suitable model may depend on the test case considered. In addition, each of these models relies on some parameters that have to be adjusted empirically, making the process of finding the right viscosity coefficient even more difficult.

A more recent approach consists in exploiting the capabilities of neural networks in pattern recognition to design data-driven tools. In short, neural networks can be described as non-linear functions with many parameters —from thousands to millions— that can be adjusted using gradient descent in order to optimize a given criterion. Examples of this approach can be found in several related topics: neural networks are used as troubled-cell indicator for high order schemes in \cite{ray2018artificial}, as classifiers of functions' regularity to control oscillations in spectral methods in \cite{schwander2021controlling}, or as predictor of the degree of reconstruction for MOOD algorithm in \cite{bourriaud2020priori}. Last but not least, in \cite{discacciati2020} the authors design an artificial viscosity for DG methods using neural networks. In all these examples, the authors use \emph{supervised learning}, where the criterion to optimize is an error between the output of the neural network and a given target. In \cite{discacciati2020} for instance, for each test case of the training dataset, the target is set to the artificial viscosity model (among a given set of options) that performs best for this specific test case. This method makes the neural network converge toward a kind of good interpolation between known models. However supervised learning in this way is not necessarily the best option, or may not even be possible in some other contexts. Indeed, an appropriate target is not always available; and when it is, as is the case in the example previously described, using it as a target will not allow the neural network to explore new designs.

In this paper we use a different way to train parameterized functions in numerical schemes, that is not bounded by these limitations. We train the neural network directly in the numerical scheme, and compare the resulting numerical solution to a target solution, instead of considering the output of the neural network itself. This approach relieves us from any prior expectation on what the output of the neural network should look like and only focuses on the result, while having the additional advantage to include effects of the neural network through many iterations of the scheme computing the gradient by automatic differentiation framework. We speak about "differentiable physic approach" \cite{thuerey2021pbdl}. This approach have been very recently used to learn discretization \cite{doi:10.1073/pnas.1814058116,dresdner2022learning}.

This approach can be see as an optimal control approach, which we design a closed-loop control of the system to fit the reference solution. A comparable approach is reinforcement learning, which is equivalent to optimal control without using the temporal transition scheme, but only examples of transitions. To our knowledge, this approach has been used for the construction of limiters in \cite{schwarz2023reinforcement} and for adjusting the weights in WENO schemes \cite{rlweno}.

One of the main ingredients of this approach is the use of deep learning frameworks (like TensorFlow or PyTorch), not only for the implementation of the neural network, but extended to the implementation of the whole numerical scheme. Indeed, the optimization of the parameters is performed with a gradient descent algorithm, which requires the computation of the gradient of the error. Since the error involves the numerical solution produced through many iterations of the parameterized numerical scheme, all the computations need to be differentiated. By implementing the numerical scheme in a deep learning framework, the computation of this highly complex gradient can be fully automated, making the optimization algorithm fairly easy to code. However, the complexity of the gradient is in itself an obstacle to the proper functionning of the algorithm, both because of its high computational cost and because of potential gradient instability. In this paper, we propose an algorithm to adress both of these limitations, and provide results of this algorithm when applied to the design of a neural network viscosity for DG schemes.

The remaining of this paper is organized as follows. In section \ref{sec:method}, we introduce a general framework for the approach we used in this work. In section \ref{sec:details}, we describe in details our implementation of this framework to the design of an artificial viscosity for DG schemes in 1D. Finally, section \ref{sec:results} gives numerical results for the advection equation, Burgers' equation and Euler's equation, with considerations on the influence of some parameters.

\section{An optimal control method for parameterized schemes}
\label{sec:method}

In this section we describe the method we use to construct an artificial viscosity. Since this method is not specific to this problem, we describe it in general terms, as a general framework to optimize parameters in a numerical scheme for partial differential equations. In our application, these parameters are the weights of a neural network designed to output the coefficient for the artificial viscosity.

\subsection{Optimization problem}

As an example, let us consider a general hyperbolic equation
\[
	\partial_t \mathbf{U}+\nabla \cdot \mathbf{F}(\mathbf{U}) = 0,
\]
with $\mathbf{U} : \mathbb{R}^d \times \mathbb{R}_+^* \rightarrow \mathbb{R}^s$ the vector of unknowns and $\mathbf{F} : \mathbb{R}^s \rightarrow \mathbb{R}^{s\times d}$ the flux of the equation.
Let us consider a given discretization of space (finite volumes, DG, WENO schemes, etc.) and time (explicit Euler, Runge-Kutta, etc.), resulting in a numerical scheme of the form
\begin{equation}
	\label{eq:scheme}
	U^{n+1} = S(U^n, \pi(U^n)) = S_{\pi}(U^n),
\end{equation}
where $U^n$ is the approximation of the solution $\mathbf{U}$ at time $t^n$, $S_{\pi}$ is an iteration of the numerical scheme on one timestep, and $\pi$ is a part of the numerical scheme that can be modified to serve as a control to the scheme. For instance, $\pi$ could be a slope limiter for a finite volumes method, a process to compute the weights of a WENO scheme, or —as in this paper— an artificial viscosity coefficient for a discontinuous Galerkin method, among many other possibilities.

In order to express the problem as an optimization problem, 
we denote by $V_{\pi}^N(U^0)$ the quantity
\begin{equation}
	V_{\pi}^N(U^0) = \mathcal{L}\big(S_{\pi}(U^0),\ S_{\pi}^2(U^0),\ \dots,\ S_{\pi}^N(U^0)\big),
\end{equation}
where $S_{\pi}^n$ corresponds to $n$ iterations of the scheme
\[
	S_{\pi}^n(U^0) = (\underbrace{S_{\pi} \circ \cdots \circ S_{\pi}}_{n\text{ times}})(U^0),
\]
and $\mathcal{L}$ is a generic cost function.
Thus $V_{\pi}^N(U^0)$ is a cost function  depending on a discrete solution made of $N$ successive iterations of the numerical scheme $S_{\pi}$. It can be for example an error committed by the scheme compared to a reference solution a penalization of the control.

Since there is usually little hope to find a control $\pi$ that minimizes $V_\pi^N(U_0)$ for all possible initial conditions $U_0$, we focus on a specific distribution of initial conditions $\P$, and consider the following optimization problem:
\begin{equation}
	\label{eq:opti-pi}
	\min_{\pi} \int V_{\pi}^N(U_0) \, d\P(U_0).
\end{equation}
In order to find a solution to this optimization problem, we choose to parameterize $\pi$ with a set of parameters $\theta$, for instance by implementing $\pi_{\theta}$ as a neural network. The optimization problem thus becomes
\begin{equation}
	\label{eq:opti-theta}
	\min_{\theta} J(\theta) = \min_{\theta} \int V_{\pi_{\theta}}^N(U_0) \, d\P(U_0).
\end{equation}
This problem is similar to 
the optimization problem solved by policy gradient methods in reinforcement learning \cite{pmlr-v32-silver14}.

\subsection{Gradient descent and back-propagation}

Our approach to solve the optimization problem \eqref{eq:opti-theta} is to use a mini-batch gradient descent algorithm, relying on automatic differentiation for the computation of the gradient. The gradient descent algorithm consists in starting from an arbitrary set of parameters, and iteratively improve it by performing updates of the form
\[
	\theta \leftarrow \theta - \eta \nabla_{\theta} J(\theta),
\]
or any other alternative, e.g with momentum, Adam, and so on. 
The \emph{mini-batch} version of the algorithm consists in replacing the gradient $\nabla_{\theta} J(\theta)$ by an approximation using a Monte-Carlo method, meaning that the integral over the distribution $\P$ of initial conditions is replaced by a sum over a sample $(U^0_1, ..., U^0_K)$ of $\P$:
\[
	\nabla_{\theta} J(\theta) \simeq \sum_{k=1}^K \nabla_{\theta} V_{\pi_{\theta}}^N(U^0_k).
\]
In principle this approximation does not prevent the convergence of the algorithm while making it much faster. In our case, if the distribution $\P$ is infinite, such an approximation is even required for the computation of the gradient to be possible.

From here, the difficulty lies in the computation of $\nabla_{\theta} V_{\pi_{\theta}}^N(U^0)$ for a given $U^0$. Figure \ref{fig:graph-V} shows the computational graph of this quantity, from which can be derived the following formulae:
\[
	\big( \nabla_{\theta} V_{\pi_{\theta}}^N \big)
	=
	\sum_{n=1}^N \big( \nabla_{\theta} S_{\pi_{\theta}}^n \big) \big( \nabla_{U^n} \mathcal{L} \big),
\]
\[
	\big( \nabla_{\theta} S_{\pi_{\theta}}^{n+1} \big)
	=
	\big( \nabla_{\theta} (S_{\pi_{\theta}} \circ S_{\pi_{\theta}}^n) \big)
	=
	\big( \nabla_{\theta} S_{\pi_{\theta}}^n \big) \big( \nabla_U S_{\pi_{\theta}} \big) + \big( \nabla_{\theta} S_{\pi_{\theta}} \big),
\]
where, for any function $g(x)$, the gradient $\nabla_x g$ refers to the transpose of its Jacobian matrix. Assuming that all these quantities are well defined, meaning that both $\mathcal{L}$ and $S_{\pi_{\theta}}$ are differentiable, we thus obtain a way to compute the gradient $\nabla_{\theta} V_{\pi_{\theta}}^N(U^0)$.

In practice, all these computations are done automatically.
Indeed, in the same way that deep learning frameworks (e.g Tensorflow, Pytorch) allow to automatically compute the gradient of a neural network w.r.t its parameters using the \emph{backpropagation} algorithm, the same frameworks can be used to compute the gradient of any parameterized numerical scheme $S_{\pi_{\theta}}$ w.r.t $\theta$, provided that $S_{\pi_{\theta}}$ is implemented using the differentiable functions of the framework. More than that, the backpropagation algorithm can be applied to any number of iterations of the numerical scheme, and even to the complete computational graph for $V_{\pi_{\theta}}^N (U^0)$ shown in Figure \ref{fig:graph-V}. In particular, this method illustrates one way these deep learning frameworks can prove useful for optimization tasks in scientific computing.

Let us conclude this section with a few observations on this approach to solve problem \eqref{eq:opti-theta}. A first advantage of this optimization algorithm lies in the fact that it does not require any reference $\pi$, since the error is not computed on the output of $\pi_{\theta}$ directly, but instead on the numerical solution that stems from $\pi_{\theta}$.
A second advantage is that the optimized function $J(\theta)$ can take into account many iterations of the numerical scheme $S_{\pi_{\theta}}$, thus including effects of the control $\pi_{\theta}$ that would go unnoticed on shorter time scales. In our application to an artificial viscosity, these long time effect would be the diffusion related to a too high viscosity, as opposed to the short term oscillations related to a too low viscosity.
Finally, note that this method contrasts with classical reinforcement learning algorithms in that this method leverages our knowledge of the transition process between two successive states, here $U^n$ and $U^{n+1}$. Indeed, classical reinforcement learning usually build implicitly (model free approaches) or explicitly (model based approaches) an approximation of this transition process by analyzing examples of transitions, whereas in our case the full knowledge of this transition process, ie of the numerical scheme, can be used to compute the gradient of $J(\theta)$ directly.

\subsection{Optimization on sub-trajectories using the reference solutions}
\label{sec:uref}

Let us mention two obstacles to the application of the method described above. The first one stems from the \emph{depth} of the computational graph when the number of iterations $N$ grows higher and higher. Indeed, as $N$ increases, the computation of $\nabla_{\theta} V_{\pi_{\theta}}^N$ becomes not only more and more expensive, but also more and more subject to gradient instability issues, similarly to very deep neural networks.
The second obstacle is that although the algorithm does not require a reference control $\pi$, it does rely on a function $\mathcal{L}$ that quantifies the error of a numerical solution, and which may not be easy to determine.

One way we have found to partially address both of these issues is to use a reference numerical scheme $S_{\text{ref}}$, accurate and robust, and the reference solutions $U_{\text{ref}}^1, ..., U_{\text{ref}}^N$ provided by it. First it helps with the measure of the error committed by $U_{\theta}$, by providing an expected result to compare with. But we can also use the reference solutions to limit the number of iterations the gradient actually goes through to a given number $m < N$, by replacing problem \eqref{eq:opti-theta} by:
\begin{equation}
	\label{eq:opti-subtraj}
	\min_{\theta} J(\theta) = \int \sum_{n=0}^{N-m} V_{\theta}^m(U_{\text{ref}}^n) \, d\P(U_0) = \min_{\theta} \int \sum_{n=0}^{N-m} V_{\theta}^m(S_{\text{ref}}^n(U_0)) \, d\P(U_0).
\end{equation}
In this formulation of the problem, instead of minimizing the error on an entire trajectory with $N$ iterations, we minimize the sum of the errors on all the sub-trajectories with $m$ iterations, starting from a point in the reference solution. This process thus limits the size of the computational graph for $J(\theta)$, while still allowing the parameterized scheme $S_{\pi_{\theta}}$ to be trained on data at times arbitrarily far from $t=0$, which would not be the case if we simply picked a small $N$. Note that for the computation of $\nabla_{\theta} J(\theta)$ with this new formulation, the sum over the sub-trajectories is also approximated by a Monte-Carlo method, similarly to the integral over the initial conditions. This approach does not allow to capture very long time effects of the control but only on medium size time sequences. For a number of applications such as the construction of limiters or viscosity this seems to be sufficient.

\subsection{Algorithm}

The whole method is described in Algorithm~\ref{alg:optimization}. Note that,  as mentioned in the previous section, both the integral over the initial conditions and the sum over the sub-trajectories of length $m$ in \eqref{eq:opti-subtraj} are approximated by a Monte-Carlo method, resulting in a kind of double mini-batch gradient descent algorithm.
Something not mentioned in this algorithm for the sake of simplicity, but very useful to track the progression of the training, is the computation of a \emph{validation loss} at the end of each epoch, consisting in the evaluation of $J(\theta)$ on a set of sub-trajectories generated at the beginning of the training.
Also note that in this algorithm, $S_{\pi_{\theta}}$ and $S_{\text{ref}}$ could actually consist of several iterations of the corresponding numerical schemes, so that the actual timestep $\Delta t$ satisfies some stability conditions. Equivalently, we could say that the error $\mathcal{L}(U^n, \cdots, U^{n+m})$ could be computed on a subset of instants, thus lowering the memory requirements for the storage of the reference solutions.

\begin{algorithm}
\caption{training algorithm}\label{alg:optimization}
\DontPrintSemicolon
Start from a random set of parameters $\theta$\;
\For{each episode}{
	Generate random initial conditions $(U^0_1, \dots, U^0_K) \sim \P$\;
	Compute reference trajectories from $U^0_k$ up to $S_{\text{ref}}^N(U^0_k)$ for all $k\in\{1, \dots, K\}$\;
	\For{each epoch}{
		Randomly select a set $I$ of indices $(k, n) \in \{1, \dots, K\}\times\{0, \dots, N\!-\!m\}$\;
		Compute sub-trajectories from $S_{\text{ref}}^n(U^0_k)$ up to $S_{\pi_{\theta}}^{m}(S_{\text{ref}}^n(U^0_k))$ for all $(k, n)\in I$\;
		Compute $J(\theta) = \sum_{(k, n)\in I} V_{\pi_{\theta}}(S_{\text{ref}}^n(U^0_k))$\;
		Update parameters $\theta$ with $\nabla J(\theta)$\;
	}
}
\end{algorithm}

\section{Design of an artificial viscosity for discontinuous Galerkin schemes}
\label{sec:details}
	
This section is dedicated to our application of the method previously described to the design of an artificial viscosity for discontinuous Galerkin schemes in one dimension.
Sections \ref{sec:visc-intro} and \ref{sec:visc-L} describe the problem and the key elements of the method, like the numerical scheme, the control $\pi$, the cost function $\mathcal{L}$.
Then sections \ref{sec:nn} to \ref{sec:scheme} give some details on the implementation.

\subsection{Discontinuous Galerkin method and artificial viscosity}
\label{sec:visc-intro}
	
In this application, we are interested in using discontinuous Galerkin (DG) schemes to solve hyperbolic equations of the form
\begin{equation}\label{eq:sys-hyp}
	\partial_t \U + \partial_x \F(\U) = 0,
\end{equation}
with $\U: \mathbb{R}_+\!\times\![x_{\min}, x_{\max}] \rightarrow \mathbb{R}^s$ the conservative variables, and $\F : \mathbb{R}^s \rightarrow \mathbb{R}^s$ the physical flux.

In order to discretize equation \eqref{eq:sys-hyp} with a discontinuous Galerkin method, we consider a spatial mesh of the interval $[x_{\min}, x_{\max}]$ made of $n_x$ cells of equal length $\Delta x$, and introduce a basis of polynomials $(\phi_1, ..., \phi_p)$ of degree at most $p-1$ on the reference interval $[-1,1]$. Assuming that the components of $\mathbf{U}$ are polynomials of degree at most $p-1$ on each cell, with no constraint of continuity at the interfaces of the cells, the $i$-th variable on the $j$-th cell can be written
\[
	\mathbf{U}_{i, j}(x, t) = \sum_{k=1}^p U_{i, j, k} (t) \phi_{k}(\hat x), \quad 1 \leq i \leq s, 1 \leq j \leq n_x,
\]
involving the change of variable to the reference interval: $\hat x = -1 + 2\ ((x-x_{\min})\!\mod\Delta x)/\Delta x \in [-1,1].$ Assuming that $\mathbf{F}(\mathbf{U})$ are also polynomials of degree at most $p-1$ on each cell, it has a similar decomposition with coefficients $(F_{i, j, k})$.
Then, integrating \eqref{eq:sys-hyp} against each of the $\phi_k$ leads to the following semi-discrete weak formulation:
\begin{align}\label{eq:semidiscrete}
\frac{dU}{dt} M + (F^{\star} - F S) = 0,
\end{align}
where $M = \left(\int_{-1}^1 \phi_k \phi_\ell\right)_{k, \ell}$, $S = \left(\int_{-1}^1 \phi_k \partial_x \phi_\ell\right)_{k, \ell}$, and $F^{\star}$ involves the estimated values at the interfaces of the cells, using a local Lax-Friedrichs flux. Here, the product $\frac{dU}{dt} M$ is to be understood as
\[
    \left(\frac{dU}{dt} M\right)_{i,j,k} = \sum_{\ell} \frac{dU_{i, j, \ell}}{dt} M_{\ell, k}.
\]
Finally, a Runge-Kutta method is used for the time integration of \eqref{eq:semidiscrete}.  We refer to \cite{hesthaven2007nodal} for more details. 

\begin{figure}
	\centering
	\includegraphics{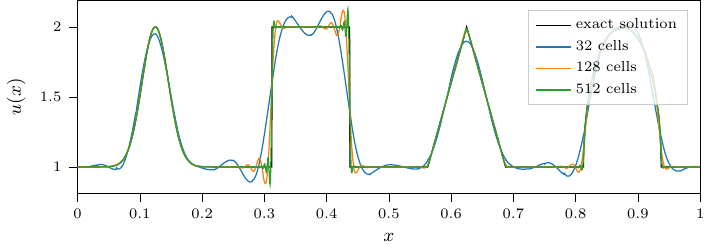}
	\caption{Example of oscillations with discontinuous Galerkin schemes. Linear advection with periodic boundary conditions, solutions after one period. All solutions were obtained using a DG scheme of order 4.}
	\label{fig:oscillations}
\end{figure}

An important benefit of discontinuous Galerkin schemes is that they can be made to converge in $O(\Delta x^p)$ for any arbitrary order $p$, by using polynomials of high enough degree ($p-1$ in one dimension).
However, when the solution exhibits strong gradients or shocks, high-order DG schemes produce oscillations as those shown in Figure \ref{fig:oscillations}, which can ruin the accuracy of the scheme and produce fatal instabilities (e.g negative pressure in the Euler equations).
For this reason, a method that is sometimes used consists in adding an artificial viscosity term to the equation to solve, in order to smooth out the solution:
\begin{equation}\label{eq:sys-visc}
	\partial_t \U + \partial_x \F(\U) = \partial_x ( \mu \partial_x \U ).
\end{equation}
The artificial viscosity above depends on a coefficient $\mu = \mu(x, t) \in \mathbb{R}$ that can locally increase or decrease the amount of smoothing, and that is expected to vanish as the length $\Delta x$ of the cells tends to zero to recover the original equation asymptotically. In practice, since the places where the viscosity is needed depend on the solution, the viscosity coefficient is taken as a function of $\U$:
\[
	\mu = \pi(\U).
\]
Denoting $\mathbf{G} = \mu \partial_x \U = \pi(\U) \partial_x \U $ and $(G_{i, j, k})$ its coefficients in the discontinuous polynomial basis, the discontinuous Galerkin scheme now reads:
\begin{align*}
	& G M = \pi(U)\, (U^{\star} - U S), \\
	& \frac{dU}{dt} M + \big( (F^{\star} - F S) - ( G^{\star}- G S ) \big) = 0.
\end{align*}
where $U^{\star}$ and $G^{\star}$ involve the estimated values at the cells interfaces using a centered numerical flux. After time discretization, still with a Runge-Kutta method, we have thus completely defined the numerical scheme $U^{n+1} = S_\pi(U^n)$.

Function $\pi$ is the one that we intend to design with the use of the method described in the previous section. As discussed in the introduction, some models for $\pi$ can already be found in the literature, and \cite{Yu:276634} compare some of them. In the result section, we compare our own viscosity to two of these models, referred to as the derivative-based (DB) and highest modal decay (MDH) models respectively, briefly described in appendix \ref{sec:viscosities}.


\subsection{Definition of the cost function}
\label{sec:visc-L}

To design the cost function $\mathcal{L}$ used to determine the control $\pi$, we compare the associated numerical solution $U^1, \cdots, U^N$ (or of a sub-trajectory $U^{n}, \cdots, U^{n+m}$) to a reference solution $U_{\text{ref}}$  by summing a local-in-time cost function $C$ over iterations:
\[
	\mathcal{L}(U^1, \dots, U^N) = \sum_{n=1}^{N} C(U^n, U_{\text{ref}}^n).
\]
For the reference solution, we use the numerical solution of a second-order MUSCL scheme on a fine grid, which ensures that the reference solution is both accurate and oscillation-free. The local-in-time cost function $C$, also concerned with both the accuracy of the solution and the presence of oscillations, is taken as a combination of three terms:
\[
	C(U^n, U_{\text{ref}}^n) =
		\omega_{\text{osc}}\, C_{\text{osc}}(U^n, U_{\text{ref}}^n)
		+ \omega_{\text{acc}}\, C_{\text{acc}}(U^n, U_{\text{ref}}^n)
		+ \omega_{\text{visc}}\, C_{\text{vis}}(U^n).
\]
For simplicity, we give below the expression of each term in the scalar case. In case of a system, we simply take the average cost.

The aim of the first term is to detect the numerical oscillations. After some testing, we have obtained interesting results with the following $W^{2,1}$ semi-norm:
\[
	C_{\text{osc}}(U^n, U_{\text{ref}}^n) = \Delta x_{\text{ref}} \sum_i \left\| D_{xx} (\Pi_{\text{ref}}(U^n))_i - D_{xx}(U_{\text{ref}}^n)_i \right\|_1,
\]
where, for any approximate quantity $\mathbf{V}$, $\mathbf{V}_i$ refers to its value in cell $i$, $\Pi_{\text{ref}}(U^n)$ is the projection of the piecewise polynomial solution of the DG scheme on the fine mesh of the reference FV scheme, and $D_{xx}(U)_i = \tfrac{1}{\Delta x_{\text{ref}}^2}(U_{i-1} - 2U_i + U_{i+1})$ a finite-difference second derivative. One can obviously use other measures of oscillations or costs penalizing positivity losses or violations of the local maximum principle.

The second term measures the accuracy of the scheme and is given by the discrete $L^1$ norm of the difference between $U^i$ and $U_{\text{ref}}^i$:
\[
	C_{\text{acc}}(U^n, U_{\text{ref}}^n) = \Delta x_{\text{ref}} \sum_i \left\| \Pi_{\text{ref}}(U^n)_i - (U_{\text{ref}}^n)_i \right\|_1,
\]
We compare the two solutions on the fine grid in order to highlight the oscillations.

Finally, since the artificial viscosity is a non-physical process, it is natural to look for the smallest viscosity which still allows to kill the oscillations. To do so, we use as the third term an $L^2$ penalisation:
\[
	C_{\text{vis}}(U^n) = \left\| \pi_{\theta}(U^n) \right\|_2^2,
\]
with the norm computed directly from the piecewise polynomial viscosity. This cost is standard in optimal control problems.

Finally, a good starting point for the weights $\omega_{\text{osc}}$, $\omega_{\text{acc}}$ and $\omega_{\text{visc}}$ could be such that all three terms contribute about as much to the overall error, but further empirical tweaking is necessary to get to the best compromise between diffusion and oscillations, as illustrated in the result section.

\subsection{Neural network viscosity function}
\label{sec:nn}

To apply Algorithm \ref{alg:optimization}, it remains to define the neural network used for the viscosity function $\pi_\theta(U)$ and how it is used in the scheme $S_\pi$.

\begin{figure}
	\centering
	\includegraphics{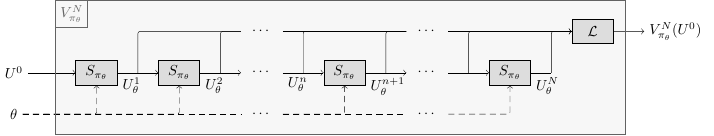}
	\caption{View of the computational graph for $V_{\pi_{\theta}}$.}
	\label{fig:graph-V}
	\vspace{3\baselineskip}
	\includegraphics{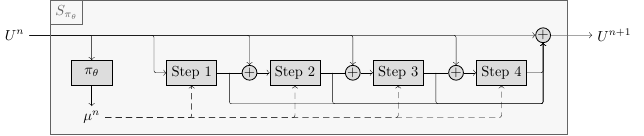}
	\caption{Computations graph for one iteration of the Runge-Kutta 4 scheme. The artificial viscosity is computed only once and used at each step.}
	\label{fig:rk4}
	\vspace{3\baselineskip}
	\includegraphics{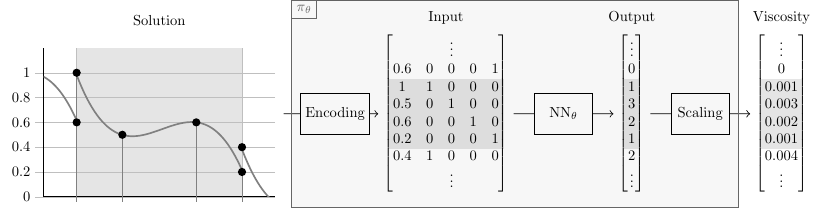}
	\caption{Computing of the viscosity with a neural network. The input variable is encoded and given to the neural network, whose output is scaled to produce the artificial viscosity.}
	\label{fig:processing}
	\vspace{3\baselineskip}
	\includegraphics{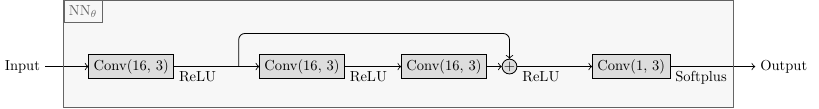}
	\caption{The architecture we use to compute the artificial viscosity: a small ResNet with only one block. Conv($w$, $k$) represents a 1D convolution with $w$ filters of size $k$.}
	\label{fig:architecture}
\end{figure}

\paragraph{Neural network architecture.} In this work we use a residual neural networks (ResNet) as introduced in \cite{he2016deep}, with adequate padding and no pooling so that the size of the output is the same as the input. This is a standard architecture for deep convolutionnal neural network. The hyper-parameters of this architecture are its \emph{depth}, i.e. the number of blocks, its \emph{width}, i.e. the number of filters per convolution, and the \emph{kernel size} of these convolutions. We got good results with a very small version of it depicted in Figure \ref{fig:architecture}: one block, width 16 and kernel size 3, for a total of about 2000 trainable parameters. We use the rectified linear unit (ReLU) activation function, except for the last layer that uses the softplus activation function. Also the last layer is initialized with kernel zero and constant bias $-3$ so that the initial output of the neural network is a constant vector with value $\softplus(-3) \simeq 0.02$. The purpose of this initialization is to start the training with a reasonable viscosity that makes the numerical scheme stable.

\paragraph{Pre-processing and post-processing}
\label{sec:processing}

The raw input for the neural network is the approximated solution $U^n$ at a given time, which comes as a tensor of values at each quadrature point of each cell. The cells are of equal length but the quadrature points are not uniformly distributed across the cells, which results in an overall non uniform discretization of the solution. Since convolutional neural networks –as the one we use– are not adapted to non uniform discretization, the input needs to be encoded in some way before being fed to the neural network. We opted for a concatenation between the value of the solution at the quadrature point and the relative position of the said quadrature point in the cell, in the form of a one-hot encoding: the first quadrature point of the cell is mapped to the vector $(1, 0, ..., 0) \in \mathbb{R}^p$, the second to $(0, 1, 0, ..., 0)\in\mathbb{R}^p$, and so on, $p$ being the number of quadrature points per cell. Figure \ref{fig:processing} gives an example of this encoding on a single variable.

Notably, the input of the neural network does not include information about the resolution of the solution, and therefore the artificial viscosity produced by the neural network is only adapted to the resolution the neural network has been trained with. In order to use the neural network with different resolution, we multiply the ouptut by a scaling factor, the same way it is done in \cite{discacciati2020}. This scaling factor $s$ is constant across each cell and involves the size of the cell $\Delta x$ as well as the jumps of the solution at its interfaces $\llbracket U \rrbracket_L$ and $\llbracket U \rrbracket_R$:
$$ s = \min \{ \Delta x, \max \{ | \llbracket U \rrbracket _L|, |\llbracket U \rrbracket_R| \} \} $$
Also, this scaling helps the artificial viscosity getting closer to zero where the solution is smooth, which prevents unnecessary diffusion.
Figure \ref{fig:processing} depicts the whole process for the computation of the viscosity on a fictive cell.

\paragraph{Integration in the numerical scheme}
\label{sec:scheme}


Since the evaluation of the neural network is relatively expensive compared to the rest of the numerical scheme, we choose to compute the artificial viscosity only once at the beginning of the timestep, as illustrated in Figure \ref{fig:rk4}. Thus, we do not update its value at each stage of the Runge-Kutta method. We found that this simplification allowed faster computing with no perceptible loss of accuracy.

\subsection{Training data}
\label{sec:ic}

In this work we try and learn from initial conditions that have a general form, expressed as partial Fourier series:
\begin{equation}
\label{eq:trainingdata}
	U^0: x \in [0, 1] \mapsto \sum_{i=0}^{20} \tfrac{a_n}{n} \cos(2\pi n x) + \tfrac{b_n}{n} \sin(2\pi n x),
\end{equation}
with coefficients $a_n$ and $b_n$ following a uniform distribution on $[-1, 1]^s$. Of course, it is possible to use other type of dataset without difficulties. For instance, for the Euler equations (see Section~\ref{sec:Euler}), we will use this kind of initialization on the primitive variables instead of the conservative ones. Positive initial conditions can be necessary for some variables: in this case, we subtract to the above functions their minima and add a small positive value $\eps = 0.1$.


As the neural network is non-local, the learned viscosity may depend on the solutions generated during the training. In particular, if the network is trained with one particular equation, it may not perform as well on another equation.  However, it would be possible to train the network directly on several equations, even if it has not been done in this work.


	\section{Numerical Results}
\label{sec:results}

In the following three sections we give numerical results for three different equations : the advection equation, Burgers' equation and Euler's system respectively.
We give some details regarding the training and the influence of some parameters in the advection case, and then simply give the results for Burgers and Euler.

In all the numerical results presented below, we use the following parameters unless stated otherwise:
\begin{itemize}
	\item DG scheme: order $p = 4$, Gauss-Lobatto quadrature points, $32$ cells on $[0, 1]$, timestep $\Delta t = 1e-5$, RK4 discretization in time,
	\item Reference FV scheme: order 2 (MUSCL), $2048$ cells on $[0, 1]$, timestep $\Delta t = 10^{-5}$, RK2 discretization in time,
	\item Entire trajectories of $N = 4096$ iterations, sub-trajectory of $m = 512$ iterations,
	\item $K = 8$ initial conditions per episode,
	\item $20$ batches of size $16$ per episode, arbitrary high number of episodes.
\end{itemize}
Values for the weights $\omega_{\text{osc}}$, $\omega_{\text{acc}}$ and $\omega_{\text{visc}}$ will be specified for each test-cases.

\subsection{Advection equation}
We will start by validating the approach on the advection equation given by
$$
\left\{\begin{array}{l}
\partial_t \rho +a\,\partial_x\rho= 0,\\
\rho(t=0,x)=\rho_0(x),
\end{array}\right.
$$
where $\rho : \mathbb{R_+}\times [0,1] \to \mathbb{R}$ is the advected density and $a \in \mathbb{R}$ is a constant velocity that we take equal to $1$. We consider periodic boudary conditions  For the initial condition $\rho_0$, we test the viscosity on a common composite function with different kinds of discontinuities:
\[
    \rho_0(x) = 1 + \left\{
        \begin{array}{ll}
            e^{-((x-0.125)/0.03)^2} & \text{if } x < 0.25, \\
            1, & \text{if } 5/16 \leq x < 7/16, \\
            1 - \left| (x - \tfrac{5}{8})\times 16 \right| &\text{if } 9/16 \leq x < 11/16, \\
            \sqrt{1 - (16x - 14)^2} &\text{if } 13/16 \leq x < 15/16, \\
            0 & \text{otherwise}.
        \end{array}
    \right..
\]
In order to notice the effects in long time of the different viscosities, we consider the solution after two periods at $t=2$. We take advantage of the simplicity of the problem to discuss the effect of the hyper parameters of the optimization problem, i.e. the weights $\omega_{\text{osc}}$, $\omega_{\text{acc}}$ and $\omega_{\text{visc}}$ involved in the cost function and the size $m$ of the sub-trajectories. 

\begin{figure}[ht!]
	\centering
	\includegraphics{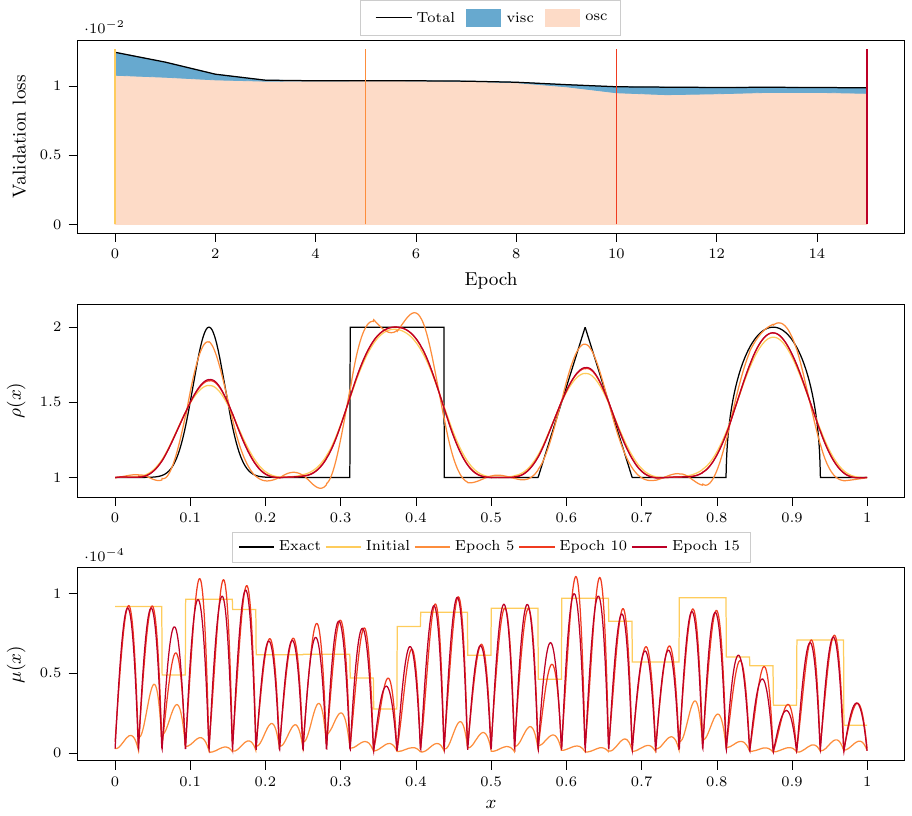}
	\caption[flushleft]{(Advection) Top: Evolution of the validation loss during training. The contribution of the different terms are shown in color. Middle and bottom: Solution (middle) and viscosity (bottom) on a test case using the neural network viscosity at different points in its training. The test case uses periodic boundary conditions and consists of two periods (ie final time $t=2$).}
	\label{fig:advection-training}
\end{figure}

For simplicity we start by training an artificial viscosity using only two of the three terms in the loss. Since the term in $C_{\text{osc}}$ on the one hand and the terms in $C_{\text{acc}}$ and $C_{\text{visc}}$ on the other seem adversarial, we consider using only $C_{\text{osc}}$ and $C_{\text{visc}}$ ($\omega_{\text{acc}}=0$), or only $C_{\text{osc}}$ and $C_{\text{acc}}$ ($\omega_{\text{visc}}=0$). Let us consider the first case. Figure \ref{fig:advection-training} shows a typical training in these conditions: at first, the neural network greatly decreases the amount of viscosity, before adding some back in order to find a better compromise between the two parts of the loss. In order to visualize the resulting viscosity and solution, the neural network is applied with a specific test case, but note that the training was done with random initial conditions as described in section \ref{sec:ic}.

\begin{figure}[ht!]
	\centering
	\includegraphics{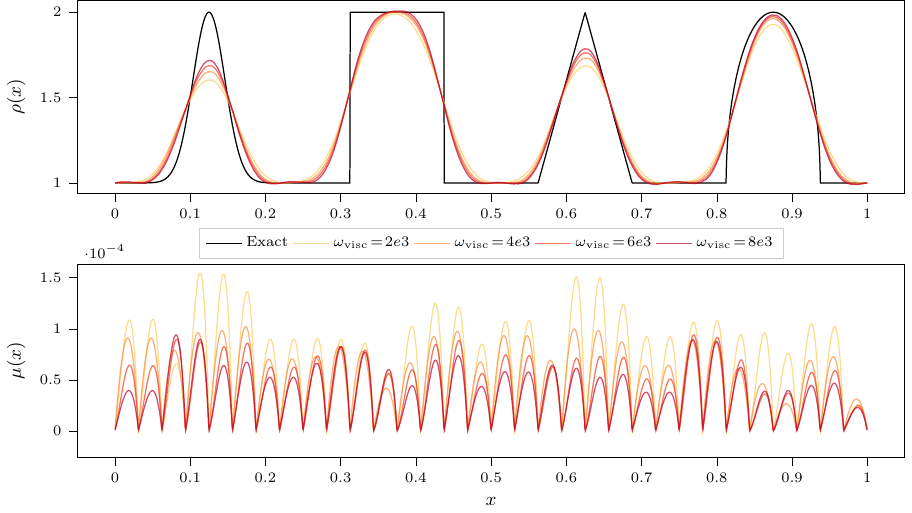}
	\caption[flushleft]{(Advection - variation $\omega_{\text{visc}}$) Solution (top) and viscosity (bottom) on a test case with neural network viscosities obtained with different weights $\omega_{\text{visc}}$, the other two weights being set to $\omega_{\text{osc}} = 10^{-5}$ and $\omega_{\text{acc}} = 0$. The test case uses periodic boundary conditions and consists of two periods (ie final time $t=2$).}
	\label{fig:advection-visc}
\end{figure}

\begin{figure}[ht!]
	\centering
	\includegraphics{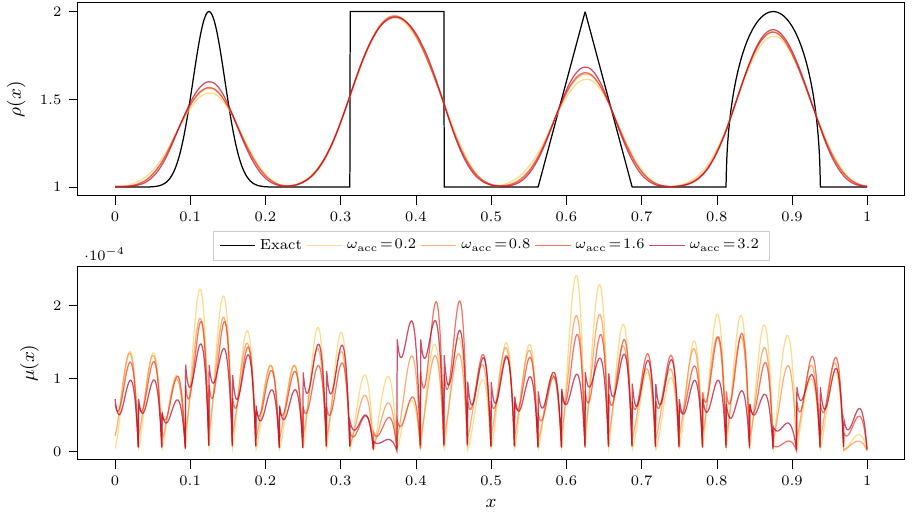}
	\caption{(Advection - variation $\omega_{\text{acc}}$) Solution (top) and viscosity (bottom) on a test case with neural network viscosities obtained with different weights $\omega_{\text{acc}}$, the other two weights being set to $\omega_{\text{osc}} = 10^{-5}$ and $\omega_{\text{visc}} = 0$. The test case uses periodic boundary conditions and consists of two periods (ie final time $t=2$).}
	\label{fig:advection-acc}
\end{figure}

An important factor in the equilibrium reached is the value chosen for the weights $\omega_{\text{osc}}$, $\omega_{\text{visc}}$ and $\omega_{\text{acc}}$. Figure \ref{fig:advection-visc} illustrates this by showing the resulting solutions and viscosities when $\omega_{\text{osc}}$ is set to $10^{-5}$ and when $\omega_{\text{visc}}$ varies ($\omega_{\text{acc}}$ still being set to zero). As expected, when the $L^2$ penalization increases, the viscosity gets smaller and smaller, which results in less diffusion but more oscillations. Indeed, as observed in Table~\ref{table:advection-visc}, both $C_{\text{osc}}$ and $L^\infty$ error increases with $\omega_{\text{visc}}$.  Figure \ref{fig:advection-acc} and Table~\ref{table:advection-acc} show what happens when it is $\omega_{\text{visc}}$ which is set to zero and $\omega_{\text{acc}}$ which varies, $\omega_{\text{osc}}$ still being set to $10^{-5}$. The results are similar, but show more diffusion overall.

\begin{table}[ht!]
	\centering
	\begin{tabular}{||c | c c c | c c c | c c||} 
		\hline
		Model & $\omega_{\text{osc}}$ & $\omega_{\text{acc}}$ & $\omega_{\text{visc}}$ & $C_{\text{osc}}$ & $C_{\text{acc}}$ & $C_{\text{visc}}$ & $L^2$ & $L^{\infty}$ \\ [0.5ex] 
		\hline\hline
		\multirow{4}{*}{DG NN} & $10^{-5}$ & 0 & $2\cdot10^3$ & 9.32e+03 &   1.10e-01 &   1.73e-10 &   5.91e-02 &   5.26e-01 \\
		 					   & $10^{-5}$ & 0 & $4\cdot10^3$ & 9.36e+03 &   9.32e-02 &   1.03e-10 &   5.95e-02 &   5.34e-01 \\
		 					   & $10^{-5}$ & 0 & $6\cdot10^3$ & 9.41e+03 &   8.35e-02 &   7.44e-11 &   5.97e-02 &   5.37e-01 \\
							   & $10^{-5}$ & 0 & $8\cdot10^3$ & 9.45e+03 &   7.57e-02 &   5.85e-11 &   5.98e-02 &   5.43e-01 \\
		\hline
	\end{tabular}
	\caption{(Advection - variation $\omega_{\text{visc}}$) Errors for each model presented in Figure \ref{fig:advection-visc}.}
	\label{table:advection-visc}
\end{table}

\begin{table}[ht!]
	\centering
	\begin{tabular}{||c | c c c | c c c | c c||} 
		\hline
		Model & $\omega_{\text{osc}}$ & $\omega_{\text{acc}}$ & $\omega_{\text{visc}}$ & $C_{\text{osc}}$ & $C_{\text{acc}}$ & $C_{\text{visc}}$ & $L^2$ & $L^{\infty}$ \\ [0.5ex] 
		\hline\hline
		\multirow{4}{*}{DG NN} & $10^{-5}$ & $0.2$ & 0 & 9.26e+03 &   1.38e-01 &   3.47e-10 &   5.86e-02 &   5.18e-01 \\
		 					   & $10^{-5}$ & $0.8$ & 0 & 9.27e+03 &   1.28e-01 &   2.87e-10 &   5.88e-02 &   5.26e-01 \\
		 					   & $10^{-5}$ & $1.6$ & 0 & 9.25e+03 &   1.28e-01 &   2.94e-10 &   5.88e-02 &   5.26e-01 \\
							   & $10^{-5}$ & $3.2$ & 0 & 9.27e+03 &   1.19e-01 &   2.55e-10 &   5.89e-02 &   5.22e-01 \\
		\hline
	\end{tabular}
	\caption{(Advection - variation $\omega_{\text{acc}}$) Errors for each model presented in Figure \ref{fig:advection-acc}.}
	\label{table:advection-acc}
\end{table}

\begin{figure}[ht!]
	\centering
	\includegraphics{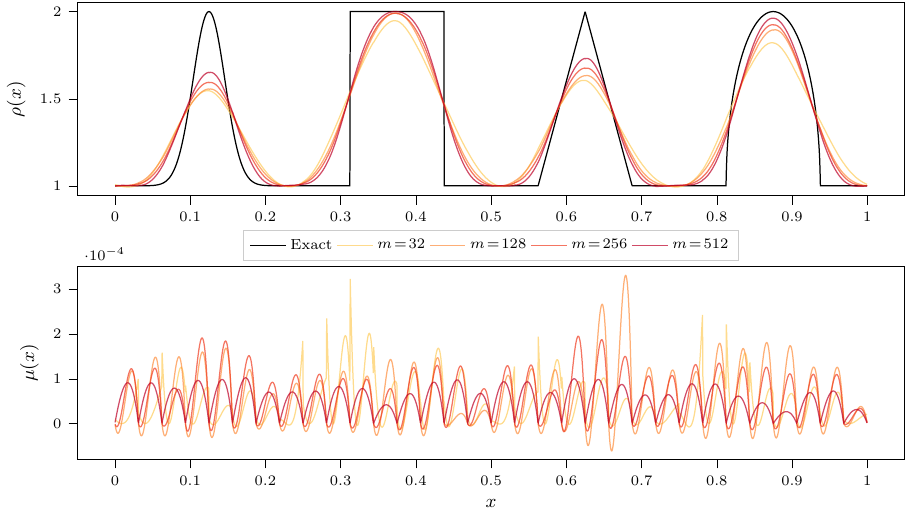}
	\caption{(Advection - variation $m$) Solution (top) and viscosity (bottom) on a test case with neural network viscosities obtained using sub-trajectories with different length $m$. The test case uses periodic boundary conditions and consists of two periods (ie final time $t=2$).}
	\label{fig:advection-watch}
\end{figure}

Another important parameter of the algorithm is the number $m$ of iterations in a sub-trajectory, on which the gradient of the loss is computed. Picking a big value for $m$ makes the computation of the gradient more expensive, but allows to include more long-term effects of the viscosity. Figure \ref{fig:advection-watch} shows some resulting viscosities with different values for $m$. It would seem that as $m$ increases, the model becomes less and less diffusive, because more diffusion appears in longer sub-trajectories and affect the gradient. In order to learn a minimal viscosity that will eliminate the oscillations it is therefore important to optimize our network on long enough trajectories so that the effect of the diffusion is clearly visible in the cost functions. Otherwise, the method learns a too large viscosity since its negative effect will not be visible enough in the cost functions. Also, a smaller $m$ do not necessarily decrease the training time even if the computation of the gradient is cheaper since the number of steps of gradient descent may be larger.

\begin{figure}[ht!]
	\centering
	\includegraphics{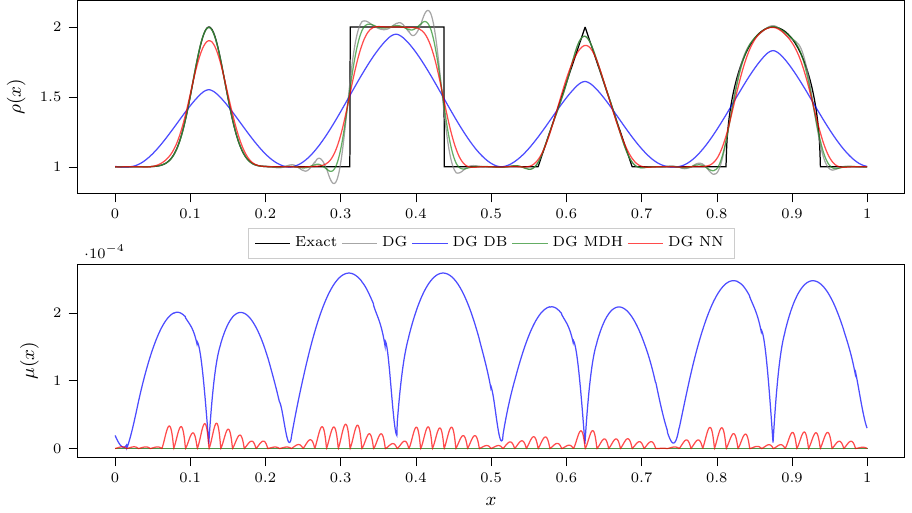}
	\caption{(Advection - comparison with DB and MDH) Solution (top) and viscosity (bottom) on a test case with 64 cells instead of the usual 32. The different viscosities used are the derivative-based (DB), the highest modal decay (MDH) and our neural network based viscosity (NN). The test case uses periodic boundary conditions and consists of two periods (ie final time $t=2$).}
	\label{fig:advection-64}
\end{figure}

\begin{figure}[ht!]
	\centering
	\includegraphics{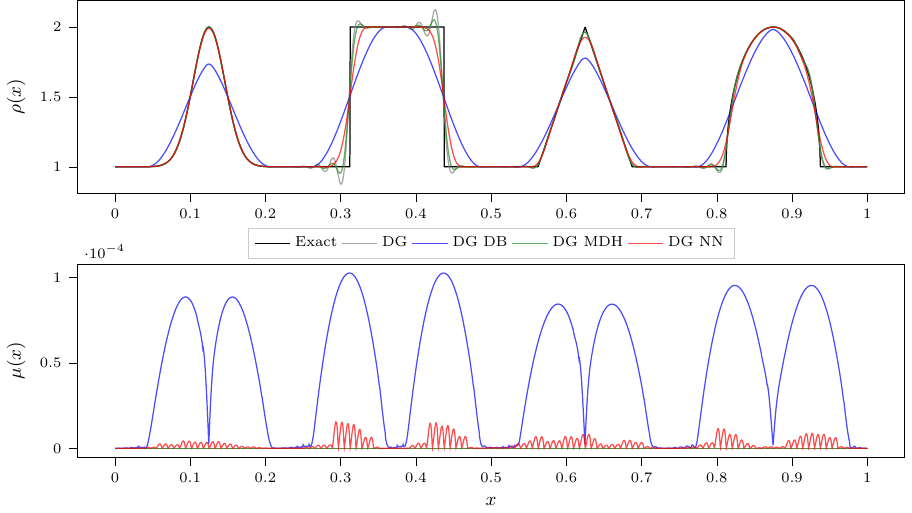}
	\caption{(Advection - comparison with DB and MDH) Solution (top) and viscosity (bottom) on a test case with 128 cells instead of the usual 32. The different viscosities used are the derivative-based (DB), the highest modal decay (MDH) and our neural network based viscosity (NN). The test case uses periodic boundary conditions and consists of two periods (ie final time $t=2$).}
	\label{fig:advection-128}
\end{figure}

\begin{figure}[ht!]
	\centering
	\includegraphics{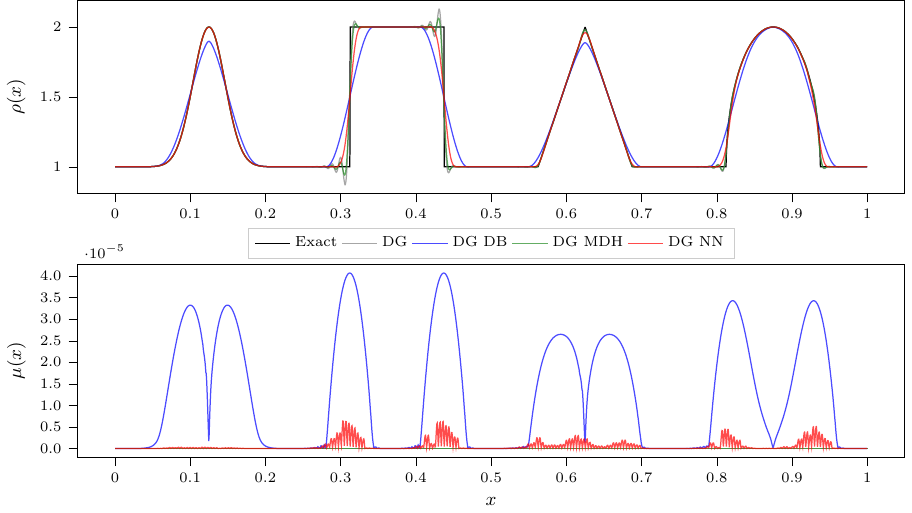}
	\caption{(Advection - comparison with DB and MDH) Solution (top) and viscosity (bottom) on a test case with 256 cells instead of the usual 32. The different viscosities used are the derivative-based (DB), the highest modal decay (MDH) and our neural network based viscosity (NN). The test case uses periodic boundary conditions and consists of two periods (ie final time $t=2$).}
	\label{fig:advection-256}
\end{figure}

Finally, we compare our viscosity to two reference viscosities: the "derivative-based" (DB) viscosity, and the "highest modal decay" (MDH) viscosity, described in appendix \ref{sec:viscosities}. Figure \ref{fig:advection-32} shows the result with $32$ cells for the discontinuous Galerkin scheme, which is the same resolution as the one used for the training of the neural network. Interestingly enough, our viscosity generalises pretty well to other resolutions, thanks to the scaling described in section \ref{sec:processing}. As an illustration, Figures \ref{fig:advection-64}, \ref{fig:advection-128} and \ref{fig:advection-256} compare the same three viscosities but used with 64, 128 and 256 cells respectively. Note that in these three examples, the model "DG NN" uses the same viscosity as in Figure \ref{fig:advection-32}, trained on $32$ cells only. The results on the different figures and given by Table~\ref{table:advection-32} show that the neural network viscosity gives the best compromise between accuracy and oscillations or between $L^2$ and $L^{\infty}$ errors. Indeed, the MDH method has a lower $L^2$ error but oscillates more and has a larger $L^{\infty}$ error. The DB approach is clearly more diffusive for this long time problem.

\begin{figure}[ht!]
	\centering
	\includegraphics{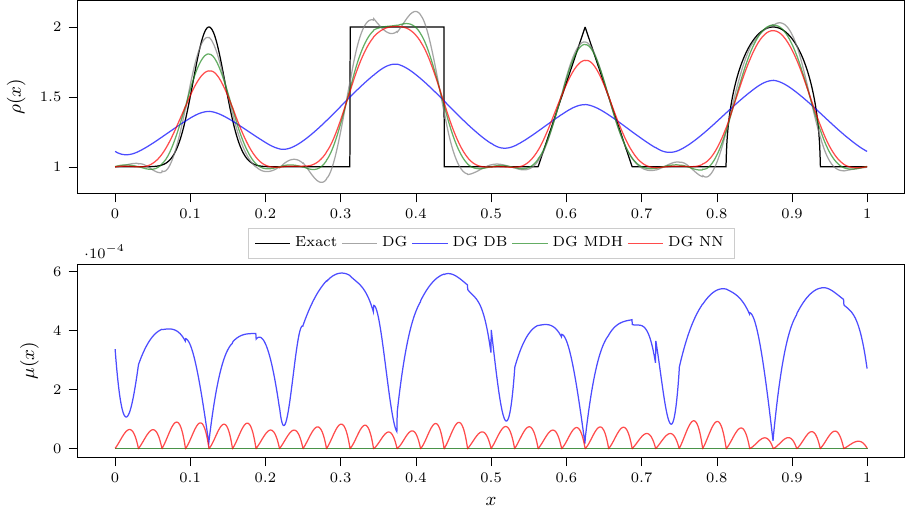}
	\caption{(Advection - comparison with DB and MDH) Solution (top) and viscosity (bottom) on a test case with different viscosities: no viscosity (DG), derivative-based viscosity (DG DB), highest modal decay viscosity (DG MDH) and our neural network based viscosity (DG NN). The test case uses periodic boundary conditions and consists of two periods (ie final time $t=2$).}
	\label{fig:advection-32}
\end{figure}

\begin{table}[ht!]
	\centering
	\begin{tabular}{|| l | r | c c c | c c||} 
		\hline
		Model & Cells & $C_{\text{osc}}$ & $C_{\text{acc}}$ & $C_{\text{visc}}$ & $L^2$ & $L^{\infty}$ \\ [0.5ex] 
		\hline\hline
		\multirow{4}{*}{DG}     & 32   &                 9.97e+03 &   5.24e-02 &   0.00e+00 &   9.03e-03 &   6.05e-01 \\
		                        & 64   &                 9.89e+03 &   2.62e-02 &   0.00e+00 &   4.71e-03 &   5.95e-01 \\
		                        & 128  &                 1.01e+04 &   1.36e-02 &   0.00e+00 &   2.52e-03 &   5.81e-01 \\
		                        & 256  &                 1.05e+04 &   7.23e-03 &   0.00e+00 &   1.37e-03 &   5.60e-01 \\
		\hline
		\multirow{4}{*}{DG DB}  & 32   &                 9.18e+03 &   2.49e-01 &   3.30e-07 &   7.82e-02 &   6.04e-01 \\
		                        & 64   &                 9.18e+03 &   1.56e-01 &   6.11e-08 &   3.94e-02 &   5.10e-01 \\
		                        & 128  &                 9.21e+03 &   8.57e-02 &   7.16e-09 &   1.86e-02 &   5.06e-01 \\
		                        & 256  &                 9.22e+03 &   4.36e-02 &   6.94e-10 &   9.13e-03 &   5.02e-01 \\
		\hline
		\multirow{4}{*}{DG MDH} & 32   &                 9.57e+03 &   5.94e-02 &   0.00e+00 &   1.22e-02 &   5.57e-01 \\
		                        & 64   &                 9.56e+03 &   2.49e-02 &   0.00e+00 &   5.37e-03 &   5.57e-01 \\
		                        & 128  &                 9.69e+03 &   1.26e-02 &   0.00e+00 &   2.75e-03 &   5.50e-01 \\
		                        & 256  &                 1.00e+04 &   6.60e-03 &   0.00e+00 &   1.45e-03 &   5.37e-01 \\
		\hline
		\multirow{4}{*}{DG NN}  & 32   &                 9.41e+03 &   8.35e-02 &   4.76e-09 &   1.78e-02 &   5.37e-01 \\
		                        & 64   &                 9.30e+03 &   4.14e-02 &   3.96e-10 &   8.48e-03 &   5.22e-01 \\
		                        & 128  &                 9.27e+03 &   2.20e-02 &   2.86e-11 &   4.98e-03 &   5.08e-01 \\
		                        & 256  &                 9.34e+03 &   1.28e-02 &   3.00e-12 &   3.11e-03 &   4.94e-01 \\
		\hline
	\end{tabular}
	\caption{(Advection - comparison with DB and MDH) Errors for each model presented in Figure \ref{fig:advection-32}.}
	\label{table:advection-32}
\end{table}

\subsection{Burgers equation}
We now consider the Burgers equation given by
$$
\left\{\begin{array}{l}
\partial_t \rho +\partial_x\left(\displaystyle \frac{\rho^2}{2}\right) = 0,\\
\rho(t=0,x)=\rho_0(x),
\end{array}\right.
$$
with $\rho : \mathbb{R}_+\times[0,1]\to \mathbb{R}$ and complemented with periodic boundary conditions. The network, the hyper-parameters and the training process are exactly the same as for the advection equation, with coefficients $\omega_{\text{acc}}=0.5$, $\omega_{\text{osc}}=10^{-5}$ and $\omega_{\text{visc}}=5$. In order to avoid any issues with non-entropic solutions the DG scheme could converge to, we only consider positive functions in the dataset. The remarks made on the hyper-parameters, the learning in the previous section on advection remains valid here. We therefore propose to give direct results comparing a learned viscosity with classical viscosities. To do this, we consider an initial condition that has not been used in the training phase of the neural network viscosity:
\begin{equation*}
\rho_0(x) = 1 + \sin(2\pi x), \quad x\in[0, 1],
\end{equation*}
with final time $t=1$.

\begin{figure}
	\centering
	\includegraphics{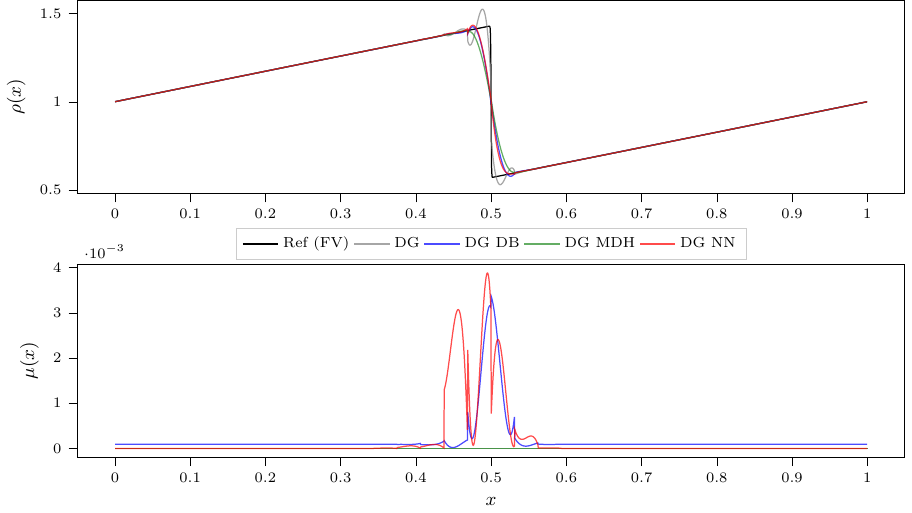}
	\caption{(Burgers - comparison with DB and MDH) Solution (top) and viscosity (bottom) of different models for Burgers equation with 32 cells}
	\label{fig:burgers-sin-32}
	\includegraphics{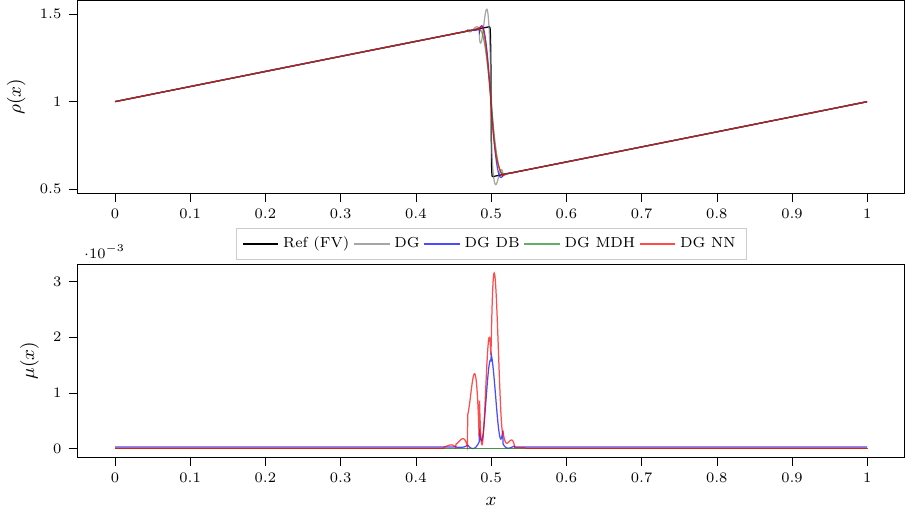}
	\caption{(Burgers - comparison with DB and MDH) Solution (top) and viscosity (bottom) of different models for Burgers equation with 64 cells}
	\label{fig:burgers-sin-64}
\end{figure}

On Figures \ref{fig:burgers-sin-32} and \ref{fig:burgers-sin-64}, we observe as before that the classical DG method without viscosity term generates large oscillations closed to the discontinuity. Contrary to the transport case, the MDH method is here the more diffusive method. Note that the MDH viscosity acts at the beginning of the simulation and then vanishes as the approximate solution becomes smooth.  The neural network and the DB methods gives very similar results with less numerical diffusion and small oscillations. Note that the neural network is slightly less oscillating at the bottom of the discontinuity. In conclusion, this test-case shows that the neural network viscosity still provides good results for such a non-linear equation with generates discontinuities.

\subsection{Euler system}
\label{sec:Euler}
Finally we present results for the Euler system:
$$
\left\{\begin{array}{l}
\partial_t \rho +\partial_x\left(\rho u\right) = 0,\\
\partial_t (\rho u) +\partial_x\left(\rho u^2 + p\right) = 0,\\
\partial_t E +\partial_x\left(Eu +pu\right) = 0,\\
\end{array}\right.
$$
where $\rho : \mathbb{R}_+\times \mathbb{R} \to \mathbb{R}$ denotes the density, $u : \mathbb{R}_+\times \mathbb{R} \to \mathbb{R}$ the velocity, $p : \mathbb{R}_+\times \mathbb{R} \to \mathbb{R}$ the pressure and $E : \mathbb{R}_+\times \mathbb{R} \to \mathbb{R}$ the energy. The system is completed with a perfect gas law, resulting in the following relation : $$E= \frac{p}{\gamma-1}+ \frac{\rho u^2}{2},$$ 
where $\gamma$ is the adiabatic constant taken equal to $1.4$ here. Once again we use the same parameters as before, the only difference being that $\pi_{\theta}$ now has three inputs, one for each conservative variable. The output is still a single viscosity coefficient $\mu(x)$, since the same viscosity is applied to each equation. The neural network is trained using the coefficients $\omega_{\text{acc}} = 0$, $\omega_{\text{osc}} = 10^{-5}$ and $\omega_{\text{visc}} = 10^3$ in the loss.  The training dataset is made using initial conditions with the three variables $\rho$, $u$, $p$ chosen according to \eqref{eq:trainingdata} with the correction to ensure the positivity of the density and the pressure.

We compare the different viscosity approaches on two classical test cases: the Sod problem and the Shu-Osher problem. In these test cases the DG scheme without viscosity is unstable and therefore is not presented. 

\begin{figure}
	\centering
	\includegraphics{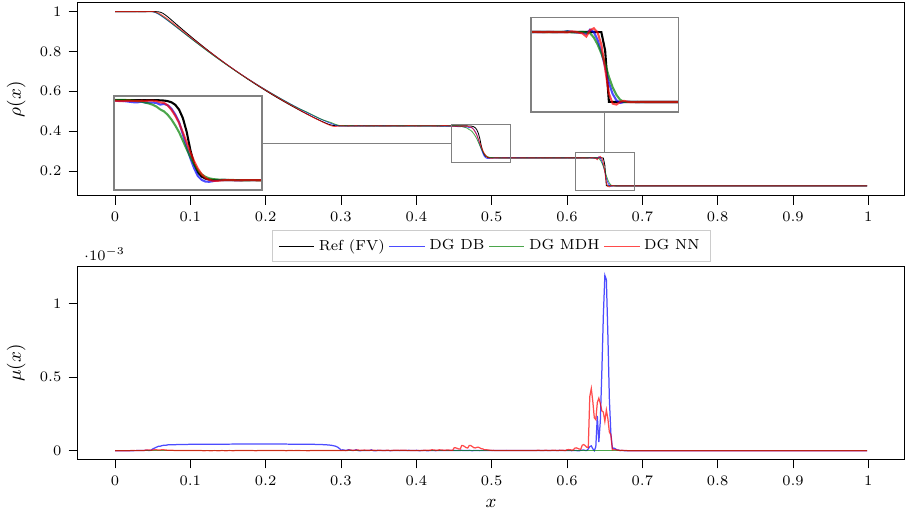}
	\caption{(Euler - Sod) Sod test case, 100 cells}
	\label{fig:euler-sod-100}
\end{figure}

\begin{figure}
	\centering
	\includegraphics{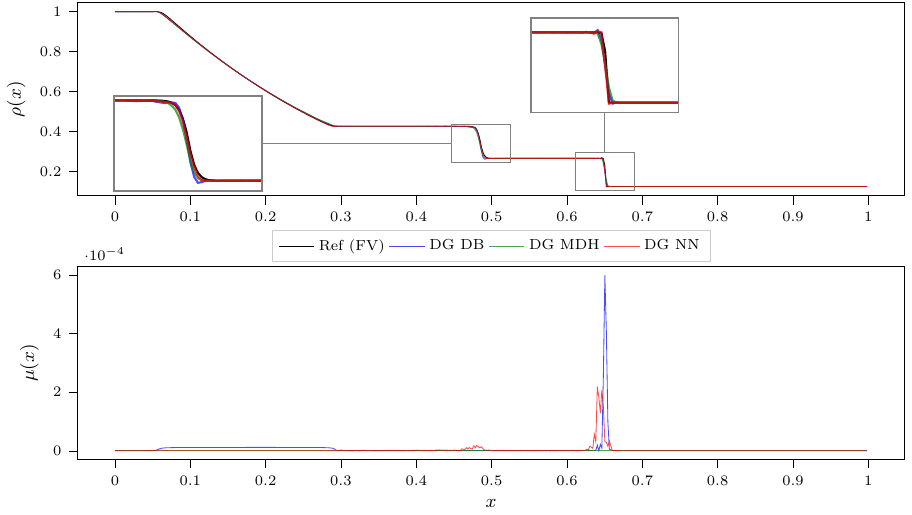}
	\caption{(Euler - Sod) Sod test case, 200 cells}
	\label{fig:euler-sod-200}
\end{figure}

The Sod test-case uses the initial condition
\begin{align*}
    (\rho_0, u_0, p_0)(x) = \left\{ \begin{array}{ll} (1, 0, 1), & \text{if } x < 0.5 \\ (0.125, 0, 0.1). & \text{otherwise} \end{array}\right.
\end{align*}
on the interval $[0,1]$ with final time $t=0.2$. We consider also Dirichlet boundary conditions.  
On Figure \ref{fig:euler-sod-100}, we compare the different schemes associated with the different viscosity models on a mesh with $100$ cells. As for the Burgers equation, the MDH viscosity provides the worst results.  This problem can be explained by the fact that the hyper-parameters of the MDH method, taken from \cite{Yu:276634}, may not be optimized to this specific test-case.  The result between the DB model and the neural network model are close. Our approach seems better in the contact wave and a little bit more oscillating on the shock. On a grid with 200 cells \ref{fig:euler-sod-200}, the neural network viscosity seems slightly more accurate for all the different components of the solution. It is confirmed by the errors presented in Table \ref{table:euler-sod}, which both $L^2$ and $L^{\infty}$ errors are smaller on the two meshes.

\begin{table}[ht!]
	\centering
	\begin{tabular}{||l | r | c c c | c c||} 
		\hline
		Model & Cells & $C_{\text{osc}}$ & $C_{\text{acc}}$ & $C_{\text{visc}}$ & $L^2$ & $L^{\infty}$ \\ [0.5ex] 
		\hline\hline
		\multirow{2}{*}{DG DB}  & 100 &  3.13e+02 &   2.47e-03 &   2.04e-08 &   3.36e-05 &   5.08e-02 \\
		                        & 200 &  3.04e+02 &   1.12e-03 &   2.48e-09 &   1.03e-05 &   3.93e-02 \\
		\hline
		\multirow{2}{*}{DG MDH} & 100 &  2.86e+02 &   2.88e-03 &   0.00e+00 &   5.27e-05 &   5.43e-02 \\
		                        & 200 &  3.00e+02 &   1.24e-03 &   0.00e+00 &   1.27e-05 &   3.90e-02 \\
		\hline
		\multirow{2}{*}{DG NN}  & 100 &  3.30e+02 &   1.36e-03 &   4.97e-09 &   1.48e-05 &   4.21e-02 \\
		                        & 200 &  2.94e+02 &   5.86e-04 &   7.09e-10 &   2.92e-06 &   3.11e-02 \\
		\hline
	\end{tabular}
	\caption{Errors for each model on the Sod test case.}
	\label{table:euler-sod}
\end{table}

\begin{figure}
	\centering
	\includegraphics{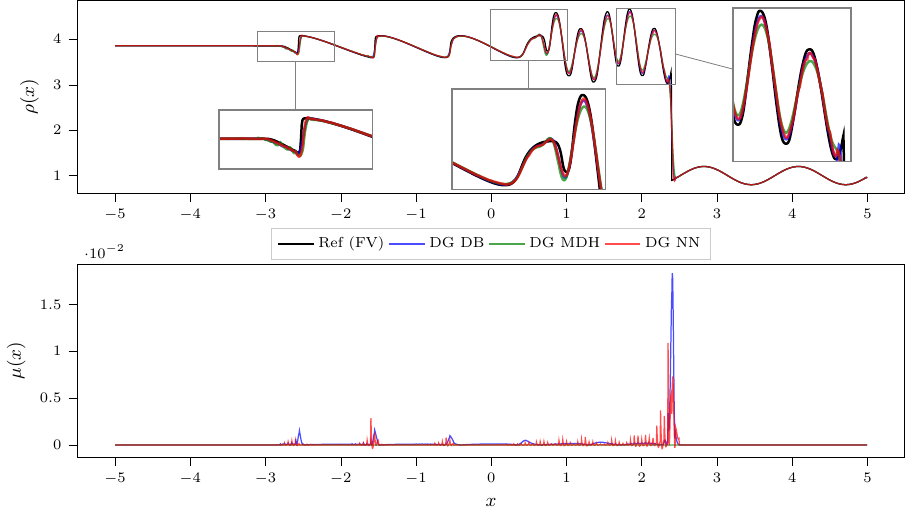}
	\caption{(Euler) Shu-Osher test case, 200 cells}
	\label{fig:euler-shuosher}
\end{figure}

The second test case is the Shu-Osher test case, whose initial condition is given by:
\begin{align*}
    (\rho_0, u_0, p_0)(x) = \left\{ \begin{array}{ll}
        (3.857143, 2.629369, 10.333333) & \text{if } x < -4 \\
        (1+0.2\sin(5x), 0, 1) & \text{otherwise}
        \end{array}\right.,
\end{align*}
on the interval $[-5,5]$ with final time $t=1.8$.
The solution is composed of several smooth oscillations and a discontinuity. As before, we compare the different approaches on a given mesh with 200 cells. The results in Figure \ref{fig:euler-shuosher} and Table \ref{table:euler-shuosher} shows that our model and the DB model gives very similar results with a slight advantage for the DB model.

\begin{figure}
	\centering
	\includegraphics{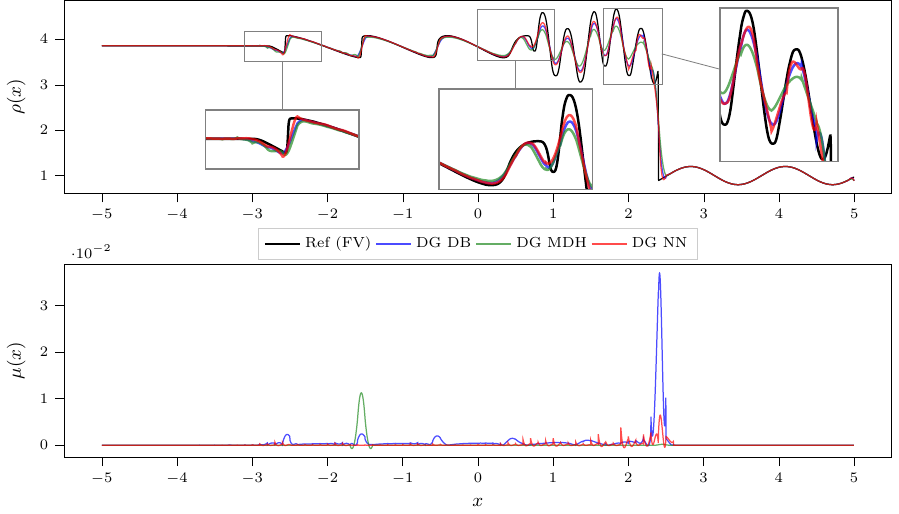}
	\caption{(Euler) Shu-Osher test case, 100 cells}
	\label{fig:euler-shuosher-100}
\end{figure}

\begin{table}[ht!]
	\centering
	\begin{tabular}{||l | r | c c c | c c||} 
		\hline
		Model & Cells & $C_{\text{osc}}$ & $C_{\text{acc}}$ & $C_{\text{visc}}$ & $L^2$ & $L^{\infty}$ \\ [0.5ex] 
		\hline\hline
		\multirow{2}{*}{DG DB}  & 100 &  2.47e+03 &   4.09e-01 &   1.73e-04 &   1.01e-01 &   1.18e+00 \\
		                        & 200 &  2.40e+03 &   1.61e-01 &   2.16e-05 &   2.92e-02 &   1.08e+00 \\
		\hline
		\multirow{2}{*}{DG MDH} & 100 &  2.37e+03 &   5.76e-01 &   1.98e-05 &   1.79e-01 &   1.31e+00 \\
		                        & 200 &  2.25e+03 &   2.53e-01 &   6.00e-13 &   5.21e-02 &   1.25e+00 \\
		\hline
		\multirow{2}{*}{DG NN}  & 100 &  2.38e+03 &   3.49e-01 &   5.46e-06 &   8.24e-02 &   1.22e+00 \\
		                        & 200 &  2.42e+03 &   1.71e-01 &   5.29e-06 &   2.95e-02 &   1.23e+00 \\
		\hline
	\end{tabular}
	\caption{(Euler) Errors for each model on the Shu-Osher test case.}
	\label{table:euler-shuosher}
\end{table}

	\section{Conclusion}

In this paper, we propose an optimal control approach to optimize a parametric numerical scheme based on its effect after several iterations. The method is a simple gradient method to optimize a given cost function, where the gradient is calculated across a large number of iterations by automatic differentiation. 
We apply it to the construction of an artificial viscosity for DG methods for one-dimensional hyperbolic equations. The numerical results on different simulations show that the obtained neural network viscosities result in equivalent or better results compared with classical artificial viscosities (Derivative Based or Highest Model Decay viscosities).

There are several possible ways to extend this work. First, non-physical oscillations have so far been detected with the semi-norm $W^{2,1}$ of the error with respect to the reference solution. Another possibility would be to design a data-driven detector of the non-physical oscillations like in \cite{beck2020neural}. 

It will also be naturally important to extend this work to 2D/3D problems. Note however that a major difficulty comes from the number of iterations taken into account in the computation of the gradient. In our one-dimensional problem, we succeed in considering up to $1000$ time steps. However, this was possible because of the coarse meshes and small networks. For two-dimensional problems, the sizes of the mesh and the network may be larger and the memory resources may be saturated. To overcome this difficulty, the method could be coupled with  a reinforcement approach \cite{rlweno} or a neural ODE method \cite{chen2018neural}, for which the gradient are computed by duality. 

Finally, the same methodology could also be applied to other problems like estimating optimal slope limiters or WENO stencils.

	\appendix

	\section{Reference artificial viscosity models}
\label{sec:viscosities}

In the result section, we compare our viscosity to two models of reference, that we briefly describe here in the context of a discontinuous Galerkin scheme of order $p$ in one dimension.

The first one is the simplest one, refered to as the derivative-based (DB) model in the comparative study \cite{Yu:276634}, and reads
\[
	\pi_{\text{DB}}(\U) = \min (\mu_{\beta}, \mu_{\max}), \quad
	\mu_{\beta} = c_{\beta} (\tfrac{\Delta x}{p-1} )^2 | \partial_x u |, \quad
	\mu_{\max} = c_{\max} \tfrac{\Delta x}{p-1} \max_{\text{cell}} | s |,
\]
where $u$ is the unique variable in the scalar case and the velocity for the Euler equation, $s$ is the local wave speed, and $c_{\beta}$ and $c_{\max}$ are empirical parameters, set to $1$ and $0.5$ respectively.

The second one is refered to as the highest modal decay (MDH) model in \cite{Yu:276634} and was first proposed in \cite{persson2006sub}. In this model, the viscosity is computed from the variable $\rho$ which refers to the unique variable in the scalar case, and to the density for the Euler equation. The MDH model relies on a modal expansion of $\rho$ in each cell,
\[
	\rho(x, t) = \sum_{k=0}^{p-1} \hat \rho_k(t) \psi_k(x), \quad \psi_k \text{ Legendre polynomials on the cell considered},
\]
and more specifically on the ratio between the norm of the highest mode and the overall norm:
\[
	r = \log_{10} \frac{\| \hat \rho_{p-1} \psi_{p-1} \|^2_{L^2}}{\| \rho \|^2_{L^2}}.
\]
The viscosity is then taken smoothly increasing with $r$ from $0$ to $\mu_{\max}$ as follows:
\[
	\pi_{\text{MDH}}(\U) = \mu_{\max} \left\{ \begin{array}{ll}
		0 & \text{if } r < r_0 - c_K \\
		\tfrac{1}{2} \left(1 + \sin \tfrac{\pi(r - r_0)}{2 c_K} \right) & \text{if } r_0 - c_K < r < r_0 + c_K \\
		1 & \text{otherwise}
		\end{array}\right.
\]
The threshold $r_0$ depends on the order $p$ as
\[
	r_0 = - \big(c_A + 4 \log_{10} (p-1) \big),
\]
and $c_A$ and $c_K$ are empirical parameters set to $2.5$ and $0.2$ respectively.
These computations give a value for the viscosity coefficient on each cell, which is interpolated by a polynomial of degree 2 that has this value in the middle of the cell, and the average value between the two cells involved at the interfaces, resulting in a continuous function.

	\bibliographystyle{alpha}
	\bibliography{sample}
\end{document}